\documentclass[12pt,twoside,a4paper,reqno]{amsart}
  
\usepackage{amsmath,amssymb} 
\usepackage{fullpage}
\usepackage{pslatex} 
\usepackage{url}

\newcommand{\Odip}[2]{\mathcal{O}_{#1}\!\Bigl(#2\Bigr)\mathchoice{\!}{}{}{}}
\newcommand{\Odi}[1]{\Odip{}{#1}}

\newcommand{\Eulerphi}{\varphi}

\begin{document}
\title{Computing the Mertens and Meissel-Mertens\\
       constants for sums over arithmetic progressions}
\author{ALESSANDRO LANGUASCO and ALESSANDRO ZACCAGNINI}
%
%
%
\begin{abstract}
We give explicit numerical values with 100 decimal digits for the
Mertens constant involved in the asymptotic formula for
$\sum\limits_{\substack{p\leq x\\ p\equiv a \bmod{q}}}1/p$
and, as a by-product, for the Meissel-Mertens constant defined as 
$\sum_{p\equiv a \bmod{q}} (\log(1-1/p)+1/p)$,
for $q \in \{3$, \dots, $100\}$ and $(q, a) = 1$.
\end{abstract}
\maketitle
\section{Introduction}
In this paper we use the technique developed in \cite{LanguascoZaccagnini2009} 
to compute the constants $M(q,a)$ involved in the following asymptotic
formula
\begin{equation}
\label{Mqa-def}
\sum\limits_{\substack{p\leq x\\ p\equiv a \bmod{q}}}
\frac{1}{p}
=
\frac{\log \log x}{\Eulerphi(q)} + M(q,a)+ \Odi{\frac{1}{\log x}},
\end{equation}
where $x\to +\infty$, 
and the so-called Meissel-Mertens constant
\[
B(q,a)
:=
\sum_{p\equiv a \bmod{q}} 
\Bigl(\log(1-\frac{1}{p})+\frac{1}{p} \Bigr),
\]
where, here and throughout the present paper,
$q \ge 3$ and $a$ are fixed integers with $(q, a) = 1$, $p$
denotes a prime number, and $\Eulerphi(q)$ is the usual Euler totient
function.
In fact we will see how to compute $M(q,a)$ with a precision of $100$
decimal digits and we will use the results in \cite{LanguascoZaccagnini2009}
to obtain the values for $B(q,a)$.

To do so we recall that the constant $C(q, a)$ studied in 
\cite{LanguascoZaccagnini2007,LanguascoZaccagnini2009} is defined implicitly by
\begin{equation}
\label{def-C}
  P(x; q, a)
  :=
  \prod_{\substack{p \le x \\ p \equiv a \bmod q}}
    \Bigl( 1 - \frac1p \Bigr)
  =
  \frac{C(q, a)}{(\log x)^{1 / \Eulerphi(q)}}
  (1 + o(1))
\end{equation}
as $x \to +\infty$.  In \cite{LanguascoZaccagnini2007} 
we proved that 
\[
  C(q, a)^{\Eulerphi(q)}
  =
  e^{-\gamma}
  \prod_p
    \Bigl( 1 - \frac1p \Bigr)^{\alpha(p; q, a)}
\]
where $\alpha(p; q, a) = \Eulerphi(q) - 1$ if $p \equiv a \bmod q$ and
$\alpha(p; q, a) = -1$ otherwise, and $\gamma$ is the Euler constant.
This enabled us to compute
their values with 100 decimal digits in \cite{LanguascoZaccagnini2009}.

Taking the logarithm of both sides in \eqref{def-C}
we get that
\[
\sum_{\substack{p \leq x \\ p\equiv a \bmod{q}}}
\log\Bigl(1-\frac{1}{p} \Bigr)
=
\log C(q,a) 
-
\frac{\log \log x}{\Eulerphi(q)}
+
o(1)
\]
as $x\to +\infty$,
and hence, adding \eqref{Mqa-def},
we obtain
\begin{equation}
\label{three-constants}
M(q,a)
=
B(q,a)
-
\log C(q,a).
\end{equation}
By \eqref{three-constants} and using the results in 
\cite{LanguascoZaccagnini2009} together
with the computation on $M(q,a)$ we will explain, we can compute
the corresponding values for $B(q,a)$ in the same range (and with the same
precision) for any $q\in \{3,\dotsc,100\}$ and $(q,a)=1$.

We recall that Finch \cite{Finch2007} has computed $M(q,a)$ and $C(q,a)$  in the case
$q \in \{ 3, 4\}$ and $(q,a)=1$.

\textbf{Acknowledgments.}
We would like to thank Robert Baillie \cite{Baillie2009}
who has driven our attention to the problem of computing $M(q,a)$.

\section{Theoretical framework}
From now on we will let $\chi$ be a Dirichlet character  $\bmod{q}$.
By the orthogonality of Dirichlet characters, a direct computation
and Theorem 428 of Hardy-Wright \cite{HardyW79} show that
\begin{equation}
\label{M-char-red}
  \Eulerphi(q)
  M(q, a)
  =
  \gamma
  +
  B
  -
 \sum_{p\mid q} 
 \frac{1}{p}
  +
  \sum_{\substack{\chi \bmod q \\ \chi \ne \chi_0}}
    \overline{\chi}(a)\sum_p \frac{\chi(p)}{p}
\end{equation}
where
\begin{equation}
\label{Meissel-Mertens-def}
 B
 :=
\sum_{p} 
\Bigl(\log(1-\frac{1}{p})+\frac{1}{p} \Bigr)
\end{equation}
is the Meissel-Mertens constant.
Moreover, using the Taylor expansion of $\log (1-x)$ and again by orthogonality, it is clear that 
\begin{equation}
\label{B-char-red}
  \Eulerphi(q)
  B(q, a)
  =
  -
  \sum_{\chi \bmod q}
    \overline{\chi}(a)
    \sum_{m \ge 2} \frac1m \sum_p \frac{\chi(p)}{p^m}
    =
  -
  \sum_{\substack{\chi \bmod q \\ \chi \ne \chi_0}}
    \overline{\chi}(a)
    \sum_{m \ge 2} \frac1m \sum_p \frac{\chi(p)}{p^m} 
    +
    B(q),  
\end{equation}
where $B(q)$,  defined as
\[
  B(q)  :  =
    -
   \sum_{m \ge 2} \frac1m \sum_{(p,q)=1} \frac{1}{p^m},
\]
represents the contribution of the principal character
$\chi_0 \bmod{q}$ and it is equal to
\[
  B(q) 
  =
  \sum_{(p,q)=1} 
  \Bigl(\log(1-\frac{1}{p})+\frac{1}{p} \Bigr)
  =
  B
-
  \sum_{p\mid q} 
\Bigl(\log(1-\frac{1}{p})+\frac{1}{p} \Bigr),
\]  
where $B$ is defined in \eqref{Meissel-Mertens-def}.
Recalling from section 2 of \cite{LanguascoZaccagnini2009} that
\begin{equation}
\label{C-char-red}
  \Eulerphi(q)
  \log C(q, a)
  =
  -\gamma
  + \log \frac{q}{\Eulerphi(q)}
  -
  \sum_{\substack{\chi \bmod q \\ \chi \ne \chi_0}}
    \overline{\chi}(a)
    \sum_{m \ge 1} \frac1m \sum_p \frac{\chi(p)}{p^m}
\end{equation}
and comparing the right hand sides of
\eqref{M-char-red}, \eqref{B-char-red} and \eqref{C-char-red},
 it is clear that it is much easier to compute $M(q,a)$ than
both $C(q,a)$ and $B(q,a)$ since in \eqref{M-char-red} no prime powers 
are involved.
Moreover, by \eqref{three-constants}, we can obtain $B(q,a)$
using $M(q,a)$ and $C(q,a)$.

Since in  \cite{LanguascoZaccagnini2009} we already computed 
several values of $C(q,a)$,  it is now sufficient to compute $M(q,a)$
for the corresponding pairs $q,a$.

To accelerate the convergence of the inner sums in 
\eqref{M-char-red}, \eqref{B-char-red} and \eqref{C-char-red},
we will consider, as we did in \cite{LanguascoZaccagnini2009}, the ``tail''
of a suitable Euler product. Letting $A$ be a fixed positive integer,
we denote the tail of the Euler product of a Dirichlet $L$-function as
\[
  L_{Aq}(\chi, s)
  =
  \prod_{p > Aq} \Bigl( 1 - \frac{\chi(p)}{p^s} \Bigr)^{-1},
\]
where $\chi \neq \chi_{0} \bmod{q}$ and $\Re(s) \ge 1$.
Now we prove that
\begin{equation}
\label{inner-sum}
  \sum_{p > A q} \frac{\chi(p)}{p^{m}}
  =
  \sum_{k \ge 1}
    \frac{\mu(k)}{k} \log(L_{A q}(\chi^k, km)),
    \end{equation}
for every integer $m \geq 1$.
We use the M\"obius inversion with a little care, since the series for
$L_{A q}(\chi, 1)$ is not absolutely convergent.
The Taylor expansion for $\log(1 - x)$ implies that
\begin{align*}
  \sum_{k \ge 2}
    \frac{\mu(k)}k \log(L_{A q}(\chi^k, k m))
  &=
  \sum_{p > A q}
    \sum_{k \ge 2}
      \sum_{n \ge 1}
        \frac{\mu(k)}{n k p^{n k m}} \chi^{n k}(p)
  =
  \sum_{p > A q}
    \sum_{\ell \ge 2}
       \frac{\chi^\ell(p)}{\ell p^{\ell m}} 
       \sum_{\substack{k \ge 2 \\ k \mid \ell}}
          \mu(k) \\
  &=
  -
  \sum_{p > A q}
    \sum_{\ell \ge 2}
       \frac{\chi^\ell(p)}{\ell p^{\ell m}}
  =
  \sum_{p > A q} \frac{\chi(p)}{p^m}
  -
  \log L_{A q}(\chi, m)
\end{align*}
since $\sum_{k \mid \ell} \mu(k) = 0$ for $\ell \ge 2$, and this proves
\eqref{inner-sum} for every $m \ge 1$.

Inserting now \eqref{inner-sum}, with $m=1$, in \eqref{M-char-red}, we have
\begin{align}
\notag
\Eulerphi(q)
&
M(q,a)
=
\gamma
  +
  B
  -
 \sum_{p\mid q} 
 \frac{1}{p}
  +
  \sum_{\substack{\chi \bmod q \\ \chi \ne \chi_0}}
    \overline{\chi}(a)
      \sum_{p \le A q} \frac{\chi(p)}{p}
    +
  \sum_{\substack{\chi \bmod q \\ \chi \ne \chi_0}}
    \overline{\chi}(a)
  \sum_{k \ge 1}
    \frac{\mu(k)}k \log(L_{A q}(\chi^k, k))
    \\
\label{M-fundamental}
&
=
\Eulerphi(q)
      \sum_{\substack{p \le A q\\ p\equiv a \bmod{q}}} \frac{1}{p}
      +
M(q)
    +
  \sum_{\substack{\chi \bmod q \\ \chi \ne \chi_0}}
    \overline{\chi}(a)
  \sum_{k \ge 1}
    \frac{\mu(k)}k \log(L_{A q}(\chi^k, k)),
\end{align}
where
\[
M(q):=
\gamma
  +
  B
  -
 \sum_{p\mid q} 
 \frac{1}{p}
-
    \sum_{\substack{p \le A q\\ (p,q)=1}} \frac{1}{p}.
\]
For $A\geq 1$, it is clear that the two sums at the 
right hand side of the previous equation collapse to
$ \sum_{p \le A q} 1/p$ but in \eqref{sum-over-a}
we will explicitly need the value of the summation over $p\mid q$ and
hence, to avoid double computations, we will use the
definition of $M(q)$ as previously stated.

For $C(q,a)$ the analogue of \eqref{M-fundamental}
is eq. (5) of \cite{LanguascoZaccagnini2009} while 
for  $B(q,a)$ it can be obtained arguing 
in a similar way.

Notice that the Riemann zeta function is never computed at $s = 1$ in
\eqref{M-fundamental}, since for $k=1$ we have
$\chi^k = \chi = \chi_0$.
To compute  the
summation over $\chi$ in \eqref{M-fundamental}
we follow the line of section 2 of \cite{LanguascoZaccagnini2009}.

This means that to evaluate \eqref{M-fundamental} using a computer program we
have to truncate the sum over $k$ and to estimate the error
we are introducing.
Let  $K > 1$ be an integer. We get
\begin{equation} 
\notag
\begin{split}
  \Eulerphi(q)
  M(q,a)
&
=
\Eulerphi(q)
      \sum_{\substack{p \le A q\\ p\equiv a \bmod{q}}} \frac{1}{p}
      +
M(q) 
      +
  \sum_{\substack{\chi \bmod q \\ \chi \ne \chi_0}}
    \overline{\chi}(a)
  \sum_{1 \le k \le K}
    \frac{\mu(k)}k \log(L_{A q}(\chi^k, k))
    \\
  &
    +
  \sum_{\substack{\chi \bmod q \\ \chi \ne \chi_0}}
    \overline{\chi}(a)
  \sum_{k > K}
    \frac{\mu(k)}k \log(L_{A q}(\chi^k, k))\\
     &=
   \widetilde{M}(q, a,A,K) + E_1(q,a,A,K),
\end{split}
\end{equation}
say.
We remark that $B$, defined as in \eqref{Meissel-Mertens-def}, 
can be easily computed up to 1000 correct digits in few seconds
by adapting  \eqref{B-char-red} to the case in which the sum in the left hand side
runs over the complete set of primes. We recall that 
Moree \cite{Moree2000}, see also the appendix by Niklasch, computed 
$B$ and many other number theoretic constants with a nice precision,
see also Gourdon-Sebah's \cite{GourdonS2001} website.
Using the Lemma in \cite{LanguascoZaccagnini2009} and the 
trivial bound for $\chi$, it is easy to see that
\[
  \left\vert E_1(q,a,A,K) \right \vert
  \leq
  \frac{2 (A q)^{1-K} (\Eulerphi(q) - 1)}{K^2 (A q - 1)}.
\]
We take this occasion to correct a typo in \cite{LanguascoZaccagnini2009}
in which, in the inequality for $E_1(q,a,A,K)$ at page 319 there,
the factor $2K$ at the denominator should be read as $K^2$.
 
In order to ensure that $\widetilde{M}(q, a,A,K)$ is a good approximation of $M(q, a)$
it is sufficient that $A q$  and $K$ are sufficiently large.
Setting $A q = 9600$ and $K = 26$ yields the desired $100$ correct
decimal digits.

Now we have to consider the error we are introducing during the
evaluation of the Dirichlet $L$-functions that appear in $\widetilde{M}(q, a,A,K)$.
This can be done exactly as in section 3 of \cite{LanguascoZaccagnini2009}
replacing $km$ there by $k$.
Let $T$ be an even integer and $N$ be a multiple of $q$.
For $\chi \neq \chi_{0} \bmod{q}$ and $k\ge 1$,
we use the Euler-MacLaurin formula
in the following form
\[
  L_{T,N}(\chi^k,k)
  =
  \sum_{r < N} \frac{\chi^k(r)}{r^{k}}
  -
  \frac1{N^{k}}
  \sum_{j = 1}^T
    \frac{(-1)^{j - 1} B_j(\chi^k)}{j!}
      \frac{k (k + 1) \cdots (k + j - 2)}{N^{j - 1}},
\]
where $B_n(\chi)$ denotes the
$\chi$-Bernoulli number which  is defined by means
of the $n$-th Bernoulli polynomial $B_n(x)$ (see Cohen 
\cite{Cohen2007b}, Definition~9.1.1), as follows
\[
  B_n(\chi)
  =
  f^{n - 1}
  \sum_{a = 0}^{f - 1} \chi(a) B_n \Bigl( \frac af \Bigr)
\]
in which $f$ is the conductor of $\chi$.

Hence the error term in evaluating the tail of the Dirichlet $L$-functions 
$L_{Aq}(\chi^k,k)$ is
\begin{align*}
  \left\vert E_2(q,a,K,N,T) \right\vert
  &\leq
  \frac{(\Eulerphi(q)-1)q^T B_T}{U(q,K,N,T)}
     \sum_{1\leq k \leq K}
       \frac{1}{k}
       \frac{k\dotsm(k+T-2)}{T!} N^{1-k-T} \\
  &=
  \frac{(\Eulerphi(q)-1)q^T B_T}{ U(q,K,N,T)T!}
     \sum_{1\leq k \leq K} (k+1)\dotsm(k+T-2) N^{1-k-T} \\
  &\leq
  \frac{(\Eulerphi(q)-1)(K+T-2)^{T-2}q^T B_T}{ U(q,K,N,T)N^{T-1}T!}
    \sum_{1\leq k \leq K} N^{-k} \\
  &\leq
  \frac{2(\Eulerphi(q)-1)(K+T-2)^{T-2}q^T B_T}{(N-1) U(q,K,N,T)N^{T-1}T!},
\end{align*}
where $B_T$ is the $T$-th Bernoulli number
and
\[
  U(q, K, N, T)
  : =
  \min_{\substack{\chi \bmod q \\ \chi \neq \chi_0}}
  \min_{1\leq k \leq K}
  \vert L_{T,N}(\chi^k,k) \vert.
\]

Collecting the previous estimates, we have that
\[
  \Bigl\vert
    M(q,a)
    -
    \frac{\widetilde{M}(q, a,A,K)}{\Eulerphi(q)}
  \Bigr\vert
  \leq
 \frac{\vert E(q,a,A,K,N,T) \vert }{\Eulerphi(q)}
\]
where $E(q,a,A,K,N,T)$ denotes
$E_1(q,a,A,K) + E_2(q,a,K,N,T)$.
 
Practical experimentations for $q \in \{3$, \dots, $100\}$ suggested us
to use different ranges for $N$ and $T$ to reach a precision of at
least $100$ decimal digits in a reasonable amount of time.
Using $A q =9600$, $K = 26$ and recalling that $q \mid N$ and $T$
is even, our choice is $N = (\lfloor 8400 / q \rfloor+1) q$ and
$T = 58$ if $q \in \{3$, \dots, $10\}$, while for $q \in \{90, \dots,100\}$ 
we have to use $N = (\lfloor 27720 / q\rfloor+1) q$ and $T = 88$.
Intermediate ranges are used for the remaining integers $q$.

The programs we used to compute the Dirichlet characters $\bmod q$ and
the values of $M(q,a)$ for $q\in \{3,\dotsc,100\}$, $1\leq a\leq q$, 
$(q,a)=1$, were written using the GP scripting language of PARI/GP
\cite{PARI2}; the C program was obtained from the GP one
using the gp2c tool. The actual computations were performed using
a double quad-core LinuX pc  for a total amount of
computing time of about 4 hours and 4 minutes.

A tiny part of the final results is collected in the tables \ref{Mfirsttable}-\ref{Bthirdtable} listed at the bottom
of this paper.
The complete set of results can be downloaded from 
\url{http://www.math.unipd.it/~languasc/Mertens-comput.html} together with the
source program in GP and the results of the verifications of the
identities \eqref{sum-over-a} and \eqref{sum-over-classes} which are
described in the section below.

Moreover, at the same web address, you will also find the values
of $B(q,a)$ computed via \eqref{three-constants} using the previous results 
on $M(q,a)$ and the ones for $C(q,a)$ in \cite{LanguascoZaccagnini2009}.
The use of \eqref{three-constants} implies some sort of ``error propagation''.
To avoid this phenomenon we recomputed some values of $C(q,a)$.
A complete report of this recomputation step can be found  at the web address previously mentioned.

Moreover, to be safer, we also directly computed $B(q,a)$ using
\eqref{B-char-red} for $q\in \{3,\dotsc,100\}$, $1\leq a\leq q$ and $(q,a)=1$. 
The needed computation time was about $3$ days, $6$ hours and a quarter.
By comparing the values of $B(q,a)$ obtained using these two different
methods, we can say that the values of $B(q,a)$ we computed are correct up to 100 decimal digits.

Finally, we also wrote a program to compute $B(q,a)$, $C(q,a)$ and $M(q,a)$
with at least 20 correct decimal digits. Comparing with \cite{LanguascoZaccagnini2009}, 
the main parameters can be chosen now in a much smaller way and so 
we were able to compute all these constants for every $3\le q \le 300$, 
$1\leq a\leq q$, $(q,a)=1$. 
In particular, the needed time on a double quad-core LinuX pc 
for the range $q\in\{3,\dotsc,200\}$ was about $5$ hours and $5$ minutes
while, for the range $q\in\{201,\dotsc,300\}$, it was about $18$ hours.
In this case we directly computed $B(q,a)$, $C(q,a)$ and $M(q,a)$ and we used
\eqref{three-constants} as a consistency check.

The whole set of these results can be downloaded at the web address previously mentioned.

\section{Verification of consistency}

The set of constants $M(q, a)$ satisfies many identities, and we
checked our results verifying that these identities hold within a very
small error.
The basic identities that we exploited are two: the first one is
\begin{equation}
\label{sum-over-a}
  \sum_{\substack{a \bmod q \\ (q, a) = 1}} M(q, a)
  =
  \gamma + B -   \sum_{p\mid q } \frac{1}{p}.
\end{equation}
This can be verified by a direct computation, taking into account the 
fact that primes dividing $q$ do not occur in any sum of the type
$\sum_{\substack{p\leq x\\ p\equiv a \bmod{q}}}
\frac{1}{p}$. 
 
The other identity is valid whenever we take two moduli $q_1$ and
$q_2$ with $q_1 \mid q_2$ and $(a, q_1) = 1$.
In this case we have
\begin{equation}
\label{sum-over-classes}
  M(q_1, a)
  =
  \sum_{\substack{j = 0 \\ (a + j q_1, q_2) = 1}}^{n - 1}
    M(q_2, a + j q_1)
    +
  \sum_{\substack{p \mid q_2 \\ p \equiv a \bmod q_1}}
     \frac1p   
\end{equation}
where $n = q_2 / q_1$.

Equation \eqref{sum-over-classes} holds also for
$B(q,a)$ with the only remark that
in the final summation 
the summand $1/p$ should be replaced by
$\log(1-1/p)+1/p)$. Concerning
\eqref{sum-over-a}, this holds for
$B(q,a)$ too  if we replace 
$\gamma - \sum_{p\mid q} 1/p$ with
$-\sum_{p\mid q} (\log(1-1/p)+1/p))$.

The proof of \eqref{sum-over-classes} 
depends on the fact that the residue class $a \bmod q_1$ is
the union of the classes $a + j q_1 \bmod q_2$, for $j \in \{0$,
\dots, $ n - 1\}$.
If $q_1$ and $q_2$ have the same set of prime factors the condition
$(a + j q_1, q_2) = 1$ is automatically satisfied, since $(a, q_1) = 1$
by our hypothesis.
On the other hand, if $q_2$ has a prime factor $p$ that $q_1$ lacks,
then there are values of $j$ such that $p \mid (a + j q_1, q_2)$ and
the corresponding value of $M(q_2, a + j q_1)$ in the right hand side
of \eqref{sum-over-classes} would be undefined.
The sum at the far right takes into account these primes.

The validity of \eqref{sum-over-a} was checked immediately at the end
of the computation of the constants $M(q, a)$, for a fixed $q$ and for
every $1 \le a \le q$ with $(q, a) = 1$ by the same program that
computed them.
These results were collected in a file and a different program checked
that \eqref{sum-over-classes} holds within a very small error by
building every possible relation of that kind for every $q_2 \in \{3$,
\dots, $100\}$ and $q_1 \mid q_2$ with $1 < q_1 < q_2$.
As in \cite{LanguascoZaccagnini2009},
the total number of identities checked is
\[
  \sum_{q = 3}^{100}
    \sum_{\substack{d \mid q \\ 1 < d < q}} \Eulerphi(d)
  =
  \sum_{q = 3}^{100} (q - 1 - \Eulerphi(q))
  =
 1907
\]
but they are not independent on one another. We did not
bother to eliminate redundancies since the total time requested for
this part of the computation is absolutely negligible.
Again as in \cite{LanguascoZaccagnini2009},
the number of independent identities is
\[
  \sum_{q = 3}^{100}
    \sum_{\substack{p \mid q \\ p < q}} \Eulerphi\Bigl( \frac qp \Bigr)
  =
  \sum_{n = 2}^{100} \pi\Bigl( \frac{100}n \Bigr) \Eulerphi(n)
  =
  1383,
\]
where $p$ denotes a prime in the sum on the left.
Please remark that in \cite{LanguascoZaccagnini2009}, page 323,
we erroneously wrote that the previous sum is equal to $1408$
which is in fact its value starting from $n=1$.

Similar checks were done also for the $20$ digits case.
Working for every $q\leq 300$ we have $12343$ independent relations 
over a total number of $17453$ ones. In this case, too,
we obtained the  desired precision (at least $20$ decimal
digits).

\vskip 1cm
\noindent
Alessandro Languasco, Dipartimento di Matematica Pura e Applicata, Universit\`a
di Padova, Via Trieste 63, 35121 Padova, Italy; languasco@math.unipd.it

\medskip
\noindent
Alessandro Zaccagnini, Dipartimento di Matematica, Universit\`a di Parma, Parco
Area delle Scienze 53/a, Campus Universitario, 43100 Parma, Italy;
alessandro.zaccagnini@unipr.it

\begin{table}
\begin{center}
\begin{tabular}{|c|c|rrr|c|} 
\hline
$q$  &  $a$  &$\ $&  \omit \hfill $M(q,a)$ \hfill  &$\ $& digits \\ \hline
3  &  1  &$\ $&   $-0.3568904795094431291196495672231858954785\dotsc$ &$\ $&104\\ 
3  &  2  &$\ $&  $0.2850543590237525795417430724985484211968\dotsc$ &$\ $&104\\ 
4  &  1  &$\ $&  $-0.2867420562261751986539451414394238573642\dotsc$ &$\ $&104\\ 
4  &  3  &$\ $&  $0.0482392690738179824093719800481197164157\dotsc$ &$\ $&104\\ 
5  &  1  &$\ $&  $-0.2088344499872589831393679436740355309848\dotsc$ &$\ $&104\\ 
5  &  2  &$\ $&  $0.3960964763519752181620428282992694487673\dotsc$ &$\ $&104\\ 
5  &  3  &$\ $&  $0.1386504040417767598465036287614177642579\dotsc$ &$\ $&104\\ 
5  &  4  &$\ $&  $-0.2644152175588502111137516747779558229888\dotsc$ &$\ $&104\\ 
\vdots  & \vdots    &$\ $& \omit \hfill \vdots \hfill  &$\ $  & \vdots    \\
9  &   1  &$\ $&   $-0.1623582321428699054929449337179721787641\dotsc$ &$\ $&104\\ 
9  &   2  &$\ $&  $0.4073663127461732280783211701614365152217\dotsc$ &$\ $&104\\ 
9  &   4  &$\ $&   $-0.1293374149143960665485300101130823600639\dotsc$ &$\ $&104\\ 
9  &   5  &$\ $&  $0.0358267016686470538569841873571831423790\dotsc$ &$\ $&104\\ 
9  &   7  &$\ $&   $-0.0651948324521771570781746233921313566504\dotsc$ &$\ $&104\\ 
9  &   8  &$\ $&   $-0.1581386553910677023935622850200712364040\dotsc$ &$\ $&104\\ 
\vdots  & \vdots    &$\ $& \omit \hfill \vdots \hfill  &$\ $  & \vdots    \\
15  &   1  &$\ $&   $-0.1506479789635675321223227319951881293922\dotsc$ &$\ $&104\\ 
15  &   2  &$\ $&  $0.3967702079831602519989220502990044120830\dotsc$ &$\ $&104\\ 
15  &   4  &$\ $&   $-0.1298987796705018718274645424839556504587\dotsc$ &$\ $&104\\ 
15  &   7  &$\ $&   $-0.0006737316311850338368792219997349633156\dotsc$ &$\ $&104\\ 
15  &   8  &$\ $&   $-0.1190129400473678821538466338276084167633\dotsc$ &$\ $&104\\ 
15  &   11  &$\ $&   $-0.0581864710236914510170452116788474015926\dotsc$ &$\ $&104\\ 
15  &   13  &$\ $&   $-0.0756699892441886913329830707443071523119\dotsc$ &$\ $&104\\ 
15  &   14  &$\ $&   $-0.1345164378883483392862871322940001725301\dotsc$ &$\ $&104\\ 
\vdots  & \vdots    &$\ $& \omit \hfill \vdots \hfill  &$\ $  & \vdots    \\
21  &   1  &$\ $&   $-0.1084483613299595805404935908928381422038\dotsc$ &$\ $&104\\ 
21  &   2  &$\ $&  $0.4250487959922326653260015663325925353478\dotsc$ &$\ $&104\\ 
21  &   4  &$\ $&   $-0.1122733018685413863141062477981428831803\dotsc$ &$\ $&104\\ 
21  &   5  &$\ $&  $0.1038169332452743207625287126777970125078\dotsc$ &$\ $&104\\ 
21  &   8  &$\ $&   $-0.0786267812146117135562454190553963721560\dotsc$ &$\ $&104\\ 
21  &   10  &$\ $&   $-0.0827607926062097370238031413241567050640\dotsc$ &$\ $&104\\ 
21  &   11  &$\ $&   $-0.0174063652116240128169859448210067632915\dotsc$ &$\ $&104\\ 
21  &   13  &$\ $&   $-0.0396915037660255136713627233956554116480\dotsc$ &$\ $&104\\ 
21  &   16  &$\ $&   $-0.0990310666212991177281328116294234061275\dotsc$ &$\ $&104\\ 
21  &   17  &$\ $&  $-0.0495316505530113215988979410810678190789\dotsc$ &$\ $&104\\ 
21  &   19  &$\ $&   $-0.0575425961745506509846081950401122043977\dotsc$ &$\ $&104\\ 
21  &   20  &$\ $&   $-0.0982465732345073585746579015543701721321\dotsc$ &$\ $&104\\ \hline
\end{tabular}
\vskip1cm
\caption{\label{Mfirsttable}
Some numerical results: the first column contains the modulus
$q$, the second the residue class $a$, the third the computed value of
$M(q,a)$ and the fourth is the number of correct decimal digits we
obtained.
The table shows the values truncated to 40 decimal digits.}
\end{center}
\end{table}

\begin{table}[htdp]
\begin{center}
\begin{tabular}{|c|c|rrr|c|}
\hline
$q$  &  $a$  &$\ $&  \omit \hfill $M(q,a)$ \hfill  &$\ $& digits \\ \hline
39  &  1 &$\ $&  $-0.0544150300747313383827161426970038945340\dotsc$ &$\ $&   104\\
39  &  2 &$\ $&  $0.4598676271292146454635405190244096879502\dotsc$ &$\ $&   104\\
39  &  4 &$\ $&  $-0.0459945989750685459192257387089186431466\dotsc$ &$\ $&   104\\
39  &  5 &$\ $&  $0.1419809783012313146832767050494466413940\dotsc$ &$\ $&   104\\
39  &  7 &$\ $&  $0.0795104580446772217182944213478570283508\dotsc$ &$\ $&   104\\
39  &  8 &$\ $&  $-0.0482711271695363892787808993829901318045\dotsc$ &$\ $&   104\\
39  &  10 &$\ $&  $-0.0625407087664201913542100549743954739736\dotsc$ &$\ $&   104\\
39  &  11 &$\ $&  $0.0351003379054991567928449805604752789995\dotsc$ &$\ $&   104\\
39  &  14 &$\ $&  $-0.0462707514061093124433385415480332066056\dotsc$ &$\ $&   104\\
39  &  16 &$\ $&  $-0.0671621258405927639177835818713861803563\dotsc$ &$\ $&   104\\
39  &  17 &$\ $&  $-0.0045382608754434839604448757618956856127\dotsc$ &$\ $&   104\\
39  &  19 &$\ $&  $-0.0078618584432586956665912002217299778289\dotsc$ &$\ $&   104\\
39  &  20 &$\ $&  $-0.0459231434298270830511730272690723063199\dotsc$ &$\ $&   104\\
39  &  22 &$\ $&  $-0.0510428073342697409080043003007927978280\dotsc$ &$\ $&   104\\
39  &  23 &$\ $&  $-0.0124631690534506350301370113904825945643\dotsc$ &$\ $&   104\\
39  &  25 &$\ $&  $-0.0581243810207459383640563019555017509681\dotsc$ &$\ $&   104\\
39  &  28 &$\ $&  $-0.0517462190695606135794914158573956691966\dotsc$ &$\ $&   104\\
39  &   29 &$\ $&  $-0.0270586216004553688028443336640559574303\dotsc$ &$\ $&   104\\
39  &   31 &$\ $&  $-0.0298803394337857678868833935975175701503\dotsc$ &$\ $&   104\\
39  &   32 &$\ $&  $-0.0473004404088889108385362179544297078717\dotsc$ &$\ $&   104\\
39  &   34 &$\ $&  $-0.0463158869206440559788110651402015355032\dotsc$ &$\ $&   104\\
39  &   35 &$\ $&  $-0.0548045111222561868022129723349733959006\dotsc$ &$\ $&   104\\
39  &   37 &$\ $&  $-0.0382400585981196219570938701692763534202\dotsc$ &$\ $&   104\\
39  &   38 &$\ $&  $-0.0652645592462251671904512528298502010370\dotsc$ &$\ $&   104\\
\hline
\end{tabular}
\vskip1cm
\caption{\label{Msecondtable}
Some numerical results: the first column contains the modulus
$q$, the second the residue class $a$, the third the computed value of
$M(q,a)$ and the fourth is the number of correct decimal digits we
obtained.
The table shows the values truncated to 40 decimal digits.}
\end{center}
\end{table}

\begin{table}[htdp]
\begin{center}
\begin{tabular}{|c|c|rrr|c|}
\hline
$q$  &  $a$  &$\ $&  \omit \hfill $M(q,a)$ \hfill  &$\ $& digits \\ \hline
84  &  1  &$\ $&  $-0.0734639142617973328342764883795225181917\dotsc$ &$\ $&  104\\
84  &  5  &$\ $&  $0.1483235495915302335618054737355898008922\dotsc$ &$\ $&  104\\
84  &  11  &$\ $&  $0.0290724024926249302145251081589204300848\dotsc$ &$\ $&  104\\
84  &  13  &$\ $&  $0.0224290728747548696540000597026220130239\dotsc$ &$\ $&  104\\
84  &  17  &$\ $&  $0.0003006811962294747858026144190460423904\dotsc$ &$\ $&  104\\
84  &  19  &$\ $&  $-0.0057630788020940875837442254844273481232\dotsc$ &$\ $&  104\\
84  &  23  &$\ $&  $-0.0138497965201530520881105236922199356889\dotsc$ &$\ $&  104\\
84  &  25  &$\ $&  $-0.0607016756429021608396776477008214514160\dotsc$ &$\ $&  104\\
84  &  29  &$\ $&  $-0.0222266883543388763218294385273009788572\dotsc$ &$\ $&  104\\
84  &  31  &$\ $&  $-0.0325213918938471846809848382364695447795\dotsc$ &$\ $&  104\\
84  &  37  &$\ $&  $-0.0442703289531425263413021626262327744388\dotsc$ &$\ $&  104\\
84  &  41  &$\ $&  $-0.0475336593998265389573863164643172518681\dotsc$ &$\ $&  104\\
84  &  43  &$\ $&  $-0.0349844470681622477062171025133156240120\dotsc$ &$\ $&  104\\
84  &  47  &$\ $&  $-0.0445066163462559127992767610577927883843\dotsc$ &$\ $&  104\\
84  &  53  &$\ $&  $-0.0464787677042489430315110529799271933763\dotsc$ &$\ $&  104\\
84  &  55  &$\ $&  $-0.0621205766407803833253627830982774246719\dotsc$ &$\ $&  104\\
84  &  59  &$\ $&  $-0.0498323317492407963847005555001138614693\dotsc$ &$\ $&  104\\
84  &  61  &$\ $&  $-0.0517795173724565634008639695556848562744\dotsc$ &$\ $&  104\\
84  &  65  &$\ $&  $-0.0611014074876142825858879099751875289632\dotsc$ &$\ $&  104\\
84  &  67  &$\ $&  $-0.0515716262256392254744286000973214317642\dotsc$ &$\ $&  104\\
84  &  71  &$\ $&  $-0.0564000928602728372344159805280953932988\dotsc$ &$\ $&  104\\
84  &  73  &$\ $&  $-0.0502394007123625523428183030876871602845\dotsc$ &$\ $&  104\\
84  &  79  &$\ $&  $-0.0547607376681565913868306490031906316886\dotsc$ &$\ $&  104\\
84  &  83  &$\ $&  $-0.0507129138346808196172715850900529202640\dotsc$ &$\ $&  104\\
\hline
\end{tabular}
\vskip1cm
\caption{\label{Mthirdtable}
Some numerical results: the first column contains the modulus
$q$, the second the residue class $a$, the third the computed value of
$M(q,a)$ and the fourth is the number of correct decimal digits we
obtained.
The table shows the values truncated to 40 decimal digits.}
\end{center}
\end{table}

\begin{table}
\begin{center}
\begin{tabular}{|c|c|rrr|c|}
\hline
$q$  &  $a$  &$\ $&  \omit \hfill $B(q,a)$ \hfill  &$\ $& digits \\ \hline
3  &   1  &$\ $&   $-0.0179374320543395898017537423354360793084\dotsc$ &$\ $&103\\
3  &   2  &$\ $&   $-0.2256492452247194384046517270072546894435\dotsc$ &$\ $&103\\
4  &   1  &$\ $&   $-0.0303152628374217668471785632748622368557\dotsc$ &$\ $&103\\
4  &   3  &$\ $&   $-0.0922560086565230005866745667406677670593\dotsc$ &$\ $&103\\
5  &   1  &$\ $&   $-0.0056989812258217866230186764730771910864\dotsc$ &$\ $&103\\
5  &   2  &$\ $&   $-0.2072541594806942995597739906452261268831\dotsc$ &$\ $&103\\
5  &   3  &$\ $&   $-0.0770818781394248684981698458665942749449\dotsc$ &$\ $&103\\
5  &   4  &$\ $&   $-0.0025398818937393664038276481789744757015\dotsc$ &$\ $&103\\
\vdots  & \vdots    &$\ $& \omit \hfill \vdots \hfill  &$\ $  & \vdots    \\
9  &  1  &$\ $& 	$-0.0020696391618847864572238027206860724807\dotsc$ &$\ $&103\\
9  &  2  &$\ $& 	$-0.1986304651091420386235919033135853530442\dotsc$ &$\ $&103\\
9  &  4  &$\ $& 	$-0.0039355992675986157162225504954464900241\dotsc$ &$\ $&103\\
9  &  5  &$\ $& 	$-0.0247398868156813518681399775437689383104\dotsc$ &$\ $&102\\
9  &  7  &$\ $& 	$-0.0119321936248561876283073891193035168035\dotsc$ &$\ $&102\\
9  &  8  &$\ $& 	$-0.0022788932998960479129198461499003980889\dotsc$ &$\ $&103\\
\vdots  & \vdots    &$\ $& \omit \hfill \vdots \hfill  &$\ $  & \vdots    \\
15  & 1 &$\ $& 	$-0.0007572379320997903470262134621931512735\dotsc$ &$\ $&103\\
15  & 2 &$\ $& 	$-0.1953264208891238586234409171756188321671\dotsc$ &$\ $&103\\
15  & 4 &$\ $& 	$-0.0016365552581033232050519597636557571823\dotsc$ &$\ $&103\\
15  & 7 &$\ $& 	$-0.0119277385915704409363330734696072947160\dotsc$ &$\ $&102\\
15  & 8 &$\ $& 	$-0.0013342030920277845401475680955985955697\dotsc$ &$\ $&103\\
15  & 11 &$\ $& $-0.0049417432937219962759924630108840398129\dotsc$ &$\ $&102\\
15  & 13 &$\ $& $-0.0036159002725660353133424956399798761365\dotsc$ &$\ $&102\\
15  & 14 &$\ $& $-0.0009033266356360431987756884153187185191\dotsc$ &$\ $&103\\
\vdots  & \vdots    &$\ $& \omit \hfill \vdots \hfill  &$\ $  & \vdots    \\
21  & 1 &$\ $& $-0.0003412956292374148069148220346920460252\dotsc$ &$\ $&103\\
21  & 2 &$\ $& $-0.1942344947334699688894003974560112287700\dotsc$ &$\ $&103\\
21  & 4 &$\ $& $-0.0002098816767539160024207752141222498501\dotsc$ &$\ $&103\\
21  & 5 &$\ $& $-0.0235093228522841911201270554154748282363\dotsc$ &$\ $&102\\
21  & 8 &$\ $& $-0.0007915893971472685099562881470728940934\dotsc$ &$\ $&102\\
21  & 10 &$\ $& $-0.0006922252022738492137317334872558725556\dotsc$ &$\ $&102\\
21  & 11 &$\ $& $-0.0046535617744727410631965497667779194062\dotsc$ &$\ $&102\\
21  & 13 &$\ $& $-0.0032481317635525831931756409952984059420\dotsc$ &$\ $&102\\
21  & 16 &$\ $& $-0.0004945238101126788996808520807420758727\dotsc$ &$\ $&102\\
21  & 17 &$\ $& $-0.0020363536140986162458739759280171063102\dotsc$ &$\ $&102\\
21  & 19 &$\ $& $-0.0016578370022937005358116763179908292912\dotsc$ &$\ $&102\\
21  & 20 &$\ $& $-0.0004239228532466525760974602939007126272\dotsc$ &$\ $&102\\
\hline
\end{tabular}
\vskip1cm
\caption{\label{Bfirsttable}
Some numerical results: the first column contains the modulus
$q$, the second the residue class $a$, the third the computed value of
$B(q,a)$ and the fourth is the number of correct decimal digits we
obtained.
The table shows the values truncated to 40 decimal digits.}
\end{center}
\end{table}

\begin{table}[htdp]
\begin{center}
\begin{tabular}{|c|c|rrr|c|}
\hline
$q$  &  $a$  &$\ $&  \omit \hfill $B(q,a)$ \hfill  &$\ $& digits \\ \hline
39  & 1 &$\ $&  $-0.0001121391210993880688819721271925997627\dotsc$ &$\ $&  102\\ 
39  & 2 &$\ $&  $-0.1934769655975371993490813769210619094240\dotsc$ &$\ $&  102\\ 
39  & 4 &$\ $&  $-0.0002995815105830464353463369791216808136\dotsc$ &$\ $&  102\\ 
39  & 5 &$\ $&  $-0.0232346770448918237834268807192861283189\dotsc$ &$\ $&  102\\ 
39  & 7 &$\ $&  $-0.0113299371144385520676547127480047428604\dotsc$ &$\ $&  102\\ 
39  & 8 &$\ $&  $-0.0002454341727558259212509705166966869744\dotsc$ &$\ $&  102\\ 
39  & 10 &$\ $&  $-0.0000459138224155805565513538531922246471\dotsc$ &$\ $&  102\\ 
39  & 11 &$\ $&  $-0.0044929480570028872659134897718135000937\dotsc$ &$\ $&  102\\ 
39  & 14 &$\ $&  $-0.0002193611514478558639221844642199225052\dotsc$ &$\ $&  102\\ 
39  & 16 &$\ $&  $-0.0000223235704255911196145171142624220894\dotsc$ &$\ $&  102\\ 
39  & 17 &$\ $&  $-0.0018338873065642236614942324023402025533\dotsc$ &$\ $&  102\\ 
39  & 19 &$\ $&  $-0.0015031092616855056695375513171955283284\dotsc$ &$\ $&  102\\ 
39  & 20 &$\ $&  $-0.0001858990560853595070479345313080934723\dotsc$ &$\ $&  102\\ 
39  & 22 &$\ $&  $-0.0001705919794542183214877183738499262132\dotsc$ &$\ $&  102\\ 
39  & 23 &$\ $&  $-0.0010537105458215435324380971901947683990\dotsc$ &$\ $&  102\\ 
39  & 25 &$\ $&  $-0.0000730986515003499343341780315856293861\dotsc$ &$\ $&  102\\ 
39  & 28 &$\ $&  $-0.0001342213885419137800239278870137043470\dotsc$ &$\ $&  102\\ 
39  & 29 &$\ $&  $-0.0006675737489094563083519829970440655189\dotsc$ &$\ $&  102\\ 
39  & 31 &$\ $&  $-0.0005852676315151717436339529174197209276\dotsc$ &$\ $&  102\\ 
39  & 32 &$\ $&  $-0.0001433829254050606314390524017643892945\dotsc$ &$\ $&  102\\ 
39  & 34 &$\ $&  $-0.0001404190990760282183009755731128149793\dotsc$ &$\ $&  102\\ 
39  & 35 &$\ $&  $-0.0000707291474709482991989950130887911099\dotsc$ &$\ $&  102\\ 
39  & 37 &$\ $&  $-0.0004011981531447411395316606101222440231\dotsc$ &$\ $&  102\\ 
39  & 38 &$\ $&  $-0.0000246764708272542810865300784362317789\dotsc$ &$\ $&  102\\ 
\hline
\end{tabular}
\vskip1cm
\caption{\label{Bsecondtable}
Some numerical results: the first column contains the modulus
$q$, the second the residue class $a$, the third the computed value of
$B(q,a)$ and the fourth is the number of correct decimal digits we
obtained.
The table shows the values truncated to 40 decimal digits.}
\end{center}
\end{table}

\begin{table}[htdp]
\begin{center}
\begin{tabular}{|c|c|rrr|c|}
\hline
$q$  &  $a$  &$\ $&  \omit \hfill $B(q,a)$ \hfill  &$\ $& digits \\ \hline
84  & 1 &$\ $&  $-0.0000119163858637686167954725330316682793\dotsc$ &$\ $&  102\\ 
84  & 5 &$\ $&  $-0.0232403602184713008048627754438543014690\dotsc$ &$\ $&  102\\ 
84  & 11 &$\ $&  $-0.0044365093956013183002165422530420512818\dotsc$ &$\ $&  102\\ 
84  & 13 &$\ $&  $-0.0032012002462998617358456975415447057261\dotsc$ &$\ $&  102\\ 
84  & 17 &$\ $&  $-0.0018681454567877949532487758242611745268\dotsc$ &$\ $&  102\\ 
84  & 19 &$\ $&  $-0.0014999868517105941255280000567662988513\dotsc$ &$\ $&  102\\ 
84  & 23 &$\ $&  $-0.0010419325887589517242899271111822879451\dotsc$ &$\ $&  102\\ 
84  & 25 &$\ $&  $-0.0000667248398877511948366874331992186367\dotsc$ &$\ $&  102\\ 
84  & 29 &$\ $&  $-0.0006749316281200272421803236576904643492\dotsc$ &$\ $&  102\\ 
84  & 31 &$\ $&  $-0.0005598064648252633420494915697566166895\dotsc$ &$\ $&  102\\ 
84  & 37 &$\ $&  $-0.0003842227265021586749603313599748870266\dotsc$ &$\ $&  102\\ 
84  & 41 &$\ $&  $-0.0003150584838202808316242249021937896620\dotsc$ &$\ $&  102\\ 
84  & 43 &$\ $&  $-0.0003293792433736461901193495016603777458\dotsc$ &$\ $&  102\\ 
84  & 47 &$\ $&  $-0.0002689626338128903152642799716205267673\dotsc$ &$\ $&  102\\ 
84  & 53 &$\ $&  $-0.0002170523788714227629800075137358681243\dotsc$ &$\ $&  102\\ 
84  & 55 &$\ $&  $-0.0000469315172527214573299434537537002158\dotsc$ &$\ $&  102\\ 
84  & 59 &$\ $&  $-0.0001682081573108212926252001037559317834\dotsc$ &$\ $&  102\\ 
84  & 61 &$\ $&  $-0.0001578501505831064102836762612245304399\dotsc$ &$\ $&  102\\ 
84  & 65 &$\ $&  $-0.0000453815847657077478783488866523727492\dotsc$ &$\ $&  102\\ 
84  & 67 &$\ $&  $-0.0001431568368661648075840877809230312134\dotsc$ &$\ $&  102\\ 
84  & 71 &$\ $&  $-0.0001166577690272412677759644893824297442\dotsc$ &$\ $&  102\\ 
84  & 73 &$\ $&  $-0.0001324187374485858716822419174992558660\dotsc$ &$\ $&  102\\ 
84  & 79 &$\ $&  $-0.0001103010836105202247205207207671888461\dotsc$ &$\ $&  102\\ 
84  & 83 &$\ $&  $-0.0001088643694263717444732353917069229651\dotsc$ &$\ $&  102\\ 
\hline
\end{tabular}
\vskip1cm
\caption{\label{Bthirdtable}
Some numerical results: the first column contains the modulus
$q$, the second the residue class $a$, the third the computed value of
$B(q,a)$ and the fourth is the number of correct decimal digits we
obtained.
The table shows the values truncated to 40 decimal digits.}
\end{center}
\end{table}

\end{document}